\documentclass[a4paper,12pt]{article}

\usepackage[english]{babel}
\usepackage{amsfonts,amsmath,amsthm,mathrsfs,amscd}
\usepackage{dsfont}
\usepackage{cite}
\usepackage{hyperref}
\usepackage{authblk}

\newtheorem{T}{Theorem}
\newtheorem{prop}{Proposition}
\newtheorem{corol}{Corollary}

\theoremstyle{remark}
\newtheorem{remark}{Remark}

\newcommand*{\cm}[1]{\mathscr{#1}}

\providecommand{\keywords}[1]
{
  \small	
  \textit{Key Words:~} #1
}

\providecommand{\MSC}[1]
{
  \small	
  \textit{Mathematics Subject Classification 2020:~} #1
}

\begin{document}

\title{Gibbs scheme in the theory of random fields}

\author{L.A. Khachatryan, B.S. Nahapetian}

\affil{\small Institute of Mathematics NAS RA, \href{mailto: linda@instmath.sci.am}{linda@instmath.sci.am}, \href{mailto: nahapet@instmath.sci.am}{nahapet@instmath.sci.am}}

\date{}

\maketitle

\begin{abstract}
The purpose of this work is to expand and clarify the concept of the class of Gibbs random fields and give its structure the form accepted in the theory of random processes. This is possible thanks to the proposed purely probabilistic definition of the Gibbs random field without referring to any physical notion. However, we do not oppose to each other the physical and probabilistic points of view on mathematical statistical physics; on the contrary, we show their natural compatibility within the framework of the suggested Gibbs scheme. The outlines of the corresponding theory are presented. In this theory, the DLR--definition is not used. At the same time, the related existence theorem provides one of adequate ways to construct Gibbs random fields (in the sense of probabilistic definition). The results of mathematical statistical physics are embedded in the theory of Gibbs random fields under development as its important (if not the most important) part. We describe in general terms the range of those problems that are more natural to solve using our definition. For example, the problem of the Gibbsian description of known classes of random fields as well as questions of validity of limit theorems of probability theory.
\end{abstract}

\keywords{Gibbs random field, conditional distribution, Gibbs form of conditional distribution, transition energy, Hamiltonian, Markov random field}

\MSC{82B03, 60G60, 60E05, 60J99}

\section*{Introduction}

The notion of a Gibbs random field with a given interaction potential was introduced by R. Dobrushin~\cite{Dobr-Gibbs} and later independently by O. Lanford and D. Ruelle~\cite{LR}. This notion (from now on referred to as the DLR--definition) naturally generalizes the classical Gibbs (Boltzmann--Gibbs) distribution, which establishes the energy-probability relationship for physical systems in finite volumes, to the case of infinite systems.

The DLR--definition was the starting point from which, in the second half of the last century, an impressive process of further mathematization of statistical physics began. The result was a broad and deep theory often called mathematical statistical physics. Among many crucial achievements of this theory, we note, for example, proposed by Ruelle~\cite{Ruelle} the description of mathematical structures of classical statistical physics underlying the thermodynamic formalism, established by Dobrushin~\cite{DUniq} the effective criterion for the presence of phase transitions in physical systems (the non-uniqueness of the Gibbs random field), as well as the Pirogov--Sinai theory~\cite{Sinai} on phase transitions in contour lattice models.

Mathematical statistical physics can be considered from two points of view --- physical and probabilistic. In the course of its creation, the DLR--definition of the Gibbs random field, using the concept of potential, was taken as a basis, which predetermined the dominance of the physical point of view in this process. As a result, mathematical statistical physics (despite the rigour of its results) acquired the structure inherent in physical theories.

As for the probabilistic point of view, it was not properly developed. The reason is clear: at that time, a purely probabilistic definition of the Gibbs random field, without refereing to any physical notion, was not formulated. It should be noted that in the absence of such a definition, many problems arising in mathematical statistical physics (being, in fact, probabilistic) forcedly were solved on the basis of the DLR--definition. This often led to rather cumbersome constructions. A typical example is the problem of the Gibbsianness of a random field.

By the probabilistic definition of a Gibbs random field we mean a definition that is based solely on its intrinsic properties, namely, on the properties of its finite-dimensional distributions. Such purely probabilistic definition of a Gibbs random field was proposed in~\cite{DN09} (see also~\cite{DN19}). In our paper, following the approach developed in~\cite{DN09, DN19}, we present an improved version of the purely probabilistic definition of a Gibbs random field (in what follows, probabilistic definition or simply $\Pr$--definition), which is more suitable for constructing the corresponding theory on its basis. Both definitions are given in terms of one-point conditional distributions. By the martingale convergence theorem, one-point finite-conditional distributions of any positive random field converge for almost all boundary conditions. According to our definition, a positive random field is Gibbsian if this convergence is uniform.

The fundamental difference between $\Pr$-- and DLR--definitions is the following. In the framework of $\Pr$--definition, we assume the random field to be given and require uniform continuity of its conditional distributions with respect to boundary conditions. In the framework of DLR--definition, an inverse problem is solved: the Gibbs specification is considered given, and the random field, the conditional distribution of which coincides almost everywhere with this specification, is called Gibbsian. It is clear that, like most inverse problems, this problem can have more than one solution, which is interpreted as a phase transition in the corresponding physical model. It should be noted that all random fields that are Gibbsian in the sense of the DLR--definition are Gibbsian in the sense of the $\Pr$--definition, too.

In this paper, based on the $\Pr$--definition, we show how to construct a theory of Gibbs random fields according to a scheme typical for the theory of random processes. Such scheme has the following structure: first, a definition of the random process is given, then its basic properties are established, and statements of a general nature are proved. After that, a fundamentally important representation theorem is given, which (being an existence theorem) solves the problem of construction of a random process on the basis of well-studied or simple non-random objects. These are, for example, stochastic matrices in the theory of Markov chains, spectral densities in the theory of stationary random processes, and covariance matrices in the theory of Gaussian processes. In the theory of Gibbs random fields under development such non-random objects are transition energy fields, with the help of which Gibbs specifications are constructed (see~\cite{DN19}). The existence theorem of Gibbs random fields with a given transition energy field appears as a representation theorem. One of the main ways to construct a transition energy field is based on the use of potential, and thereby, the results of mathematical statistical physics are naturally included in the proposed theory of Gibbs random fields, representing its important (if not the most important) part. Here we see the natural compatibility of the physical and probabilistic points of view on mathematical statistical physics.

It is well known that Gibbs distribution has a lot of applications outside mathematical statistical physics, especially in the case of random fields defined in finite volumes (see, for example,~\cite{Besag, GemGem, Winkler}). In~\cite{KhN22}, the proposed Gibbs scheme was applied for finite random fields. Particularly, it was shown that all positive finite random fields are Gibbsian.

The structure of this paper is as follows. The first section provides some general notations and definitions. The next section considers various systems of conditional distributions of a random field and the problem of reconstructing a random field from these systems. The main result of this part of the paper is that the conditional probabilities of any positive random field admit a Gibbsian representation. The third section gives a probabilistic definition of a Gibbs random field and provides a comparative analysis of the $\Pr$--definition with the DLR--definition and other existing definitions. The frames of the theory of Gibbs random fields are described, as well as the range of problems that can be solved on the basis of $\Pr$--definition. Although the work considers the case of a finite state space, the same scheme applies to the case of a continuous state space with a finite measure. The method we are developing can also be applied to random fields with non-compact state space, which will be the topic for a separate publication.

\section{Preliminaries}

Let $\mathbb{Z}^d$ be a $d$-dimensional integer lattice, i.e., a set of $d$-dimensional vectors with integer components, $d \geq 1$. Note that all the arguments in this paper remain valid if we consider an arbitrary countable set instead of $\mathbb{Z}^d$.

For $S \subset \mathbb{Z}^d$, denote by $W(S) = \{ V \subset S , 0 < \left| V \right| < \infty \}$ the set of all nonempty finite subsets of $S$, where $\left| V \right|$ is the number of points in $V$. In the case $S = \mathbb{Z}^d$, we will use the simpler notation $W$. To denote the complement of the set $S$, we will write $S^c$. When denoting one-point sets (singletons) $\{t\}$, $t \in \mathbb{Z}^d$, the brackets will be usually omitted.

A neighborhood system in $\mathbb{Z}^d$ is a collection $\partial = \{ \partial t, t \in \mathbb{Z}^d \}$ of finite subsets of the lattice $\mathbb{Z}^d$ such that $t \notin \partial t$ and $s \in \partial t$ if and only if $t \in \partial s$, $s \in \mathbb{Z}^d$.

Let each point $t \in \mathbb{Z}^d$ be associated with a set $X^t$, which is a copy of some finite set $X$, $1< \vert X \vert <\infty$, and let $\cm{B}^t = \cm{B}$ be a total $\sigma$-algebra on $X$. Denote by $X^S$ the set of all configurations on $S$, $S \subset \mathbb{Z}^d$, that is, the set $X^S = \{ ( x_t ,t \in S ) \}$, $x_t \in X$, of all functions defined on $S$ and tacking values in $X$. For $S = {\O}$, we assume that $X^{\O} = \{ \boldsymbol{{\O}} \}$ where $\boldsymbol{{\O}}$ is an empty configuration. For any disjoint $S, T \subset \mathbb{Z}^d$ and any $x \in X^S$, $y \in X^T$, denote by $xy$ the concatenation of $x$ and $y$, that is, the configuration on $S \cup T$ equal to $x$ on $S$ and to $y$ on $T$. When $T \subset S$, we denote by $x_T$ the restriction of configuration $x \in X^S$ on $T$, i.e., $x_T=(x_t, t \in T)$.

For each $S \subset \mathbb{Z}^d$, the set $X^S$ is endowed with the $\sigma$-algebra $\cm{B}^S$, which is the product of the $\sigma$-algebras $\cm{B}^t$, $t \in S$.

A \emph{random field} on $\mathbb{Z}^d$ with state space $X$ is a probability measure $P$ on $(X^{\mathbb{Z}^d}, \cm{B}^{\mathbb{Z}^d})$. By $P_S$ we will denote the restriction of $P$ on $S$, $S \subset \mathbb{Z}^d$, that is, a probability distribution on $(X^S, \cm{B}^S)$ such that
$$
P_S (A) = (P)_S (A) = P(A X^{\mathbb{Z}^d \backslash S}), \qquad A \in \cm{B}^S.
$$
For $S = {\O}$, we assume $P_{\O}(\boldsymbol{{\O}})=1$, and for $V \subset S$, we have $(P_S)_V = P_V$. By $\{P_V, V \in W\}$ we denote the system of finite-dimensional distributions of the random field $P$. Since it restores $P$, for brevity we will write $P = \{P_V, V \in W\}$.

A random field $P = \{P_V, V \in W\}$ is called \emph{positive} if for any $V \in W$, one has $P_V (x) > 0$ for all $x \in X^V$.

For any $S \subset \mathbb{Z}^d$ and any function $h: W(S) \to \mathbb{R}$, the notation $\lim \limits_{\Lambda \uparrow S} h(\Lambda) = a$ means that for any $\varepsilon > 0$ there exists $\Lambda_\varepsilon \in W(S)$ such that for any $\Lambda \in W(S)$, $\Lambda \supset \Lambda_\varepsilon$, it holds $\left| h(\Lambda) - a \right| < \varepsilon$. We will also consider limits $\lim \limits_{n \to \infty} h(\Lambda_n)$ with respect to some increasing sequence $\{\Lambda_n\}_{n \ge 1}$ of finite sets converging to $S$, that is, $\Lambda_n \subset \Lambda_{n+1} \in W(S)$, $n \ge 1$, and $\bigcup \limits_{n=1}^\infty \Lambda_n = S$. It is clear that if $\lim \limits_{\Lambda \uparrow S} h(\Lambda) = a$, then there exists an increasing sequence $\{\Lambda_n\}_{n \ge 1}$ of finite sets converging to $S$ such that $\lim \limits_{n \to \infty} h(\Lambda_n) = a$. At the same time, if $\lim \limits_{n \to \infty} h(\Lambda_n) = a$ for any increasing sequence $\{\Lambda_n\}_{n \ge 1}$, then $\lim \limits_{\Lambda \uparrow S} h(\Lambda) = a$.

Any increasing sequence $\cm{F} = \{\cm{F}_n\}_{n \ge 1}$ of $\sigma$-subalgebras of $\cm{B}^{\mathbb{Z}^d}$, that is, $\cm{F}_n \subset \cm{F}_{n+1} \subset \cm{B}^{\mathbb{Z}^d}$, $n \ge 1$, is called a filtration on $(X^{\mathbb{Z}^d}, \cm{B}^{\mathbb{Z}^d})$. In particular, $\cm{F} = \{\cm{B}^{\Lambda_n}\}_{n \ge 1}$ forms a filtration corresponding to increasing sequence $\{\Lambda_n\}_{n \ge 1}$ of finite sets converging to $\mathbb{Z}^d$. In what follows, we will consider only such filtrations, usually without explicitly specifying the corresponding increasing sequence.

A function $f: X^S \to \mathbb{R}$, $S \subset \mathbb{Z}^d$, is called \emph{local} if there exists  $\Lambda \in W(S)$ such that for $x \in X^S$, the value $f(x)$ depends only on the restriction $x_\Lambda$ of configuration $x$ on $\Lambda$, that is, $f(uy)=f(uz)$ for any $u \in X^\Lambda$ and $y,z \in X^{S \backslash \Lambda}$.

A function $f: X^S \to \mathbb{R}$, $S \subset \mathbb{Z}^d$, is called \emph{quasilocal} if it is the uniformly convergent limit of some sequence of local functions. Equivalently, $f$ is quasilocal if
$$
\lim_{\Lambda \uparrow S} \sup_{x,y \in X^S \,:\, x_{\Lambda} = y_{\Lambda}} \bigl| f(x) - f(y) \bigr| = 0.
$$
We will say that $f$ is quasilocal on $\cm{Y} \subset X^S$ if the relation above holds when supremum is taken over configurations $x,y \in\cm{Y}$ coinciding on $\Lambda$.

\section{Conditional distribution of a random field and its Gibbs form}
\label{sect-cond}

The results presented in this section are of a general nature. The most important one is that positive conditional probabilities of any random field can be represented in a Gibbs form.

\subsection{Conditional distribution}

In this paper, we will consider positive random fields only.

Let $P = \{P_V, V \in W\}$ be a random field and let
\begin{equation}
\label{P-cond-dens-fin}
g_V^z (x) = \dfrac{P_{V \cup \Lambda}(xz)}{P_\Lambda(z)}, \qquad x \in X^V,
\end{equation}
be its conditional probability on $X^V$ under condition $z \in X^\Lambda$, $\Lambda \in W(V^c)$, $V \in W$. The distributions~\eqref{P-cond-dens-fin} will be called \emph{finite-conditional distributions} of the random field $P$.

For any filtration $\cm{F} = \{\cm{B}^{\Lambda_n}\}_{n \ge 1}$, by virtue of the martingale convergence theorem, for any $V \in W$, $x \in X^V$ and $P$-a.e. $\bar x \in X^{V^c}$, there exist limits
\begin{equation}
\label{P-cond-dens-inf}
g_V^{\bar x} (x) = \lim \limits_{n \to \infty} g_V^{\bar x_{\Lambda_n \backslash V}}(x).
\end{equation}
Following the terminology accepted in mathematical statistical physics, any configuration $\bar x \in X^{V^c}$ will be called a boundary condition.

For the filtration $\cm{F}$ and any $V \in W$, denote by $\cm{X}_V(\cm{F})$ the set of all boundary conditions $\bar x \in X^{V^c}$ for which the limits in~\eqref{P-cond-dens-inf} exist for all $x \in X^V$. We will call the elements of $\cm{X}_V(\cm{F})$, $V \in W$, \emph{admissible boundary conditions} for the random field $P$ with respect to filtration $\cm{F}$. It is clear that for each $V \in W$, $P(\cm{X}_V(\cm{F}))=1$.

The distributions defined by~\eqref{P-cond-dens-inf} will be called \emph{conditional distributions} and the set $G(P, \cm{F}) = \{g_V^{\bar x}, \bar x \in \cm{X}_V(\cm{F}), V \in W\}$ will be called the \emph{system of conditional distributions of the random field~$P$ with respect to filtration $\cm{F}$}.

The following statement establishes a connection between the elements of the system of conditional distributions.

\begin{prop}
\label{admiss}
Let $G(P, \cm{F})$ be the system of conditional distributions of random field $P$ with respect to filtration $\cm{F}$. Then for any $I \subset V \in W$, $\bar x \in \cm{X}_V(\cm{F})$ and $y \in X^{V \backslash I}$, one has $\bar x y \in \cm{X}_I(\cm{F})$, and for all $x,u \in X^I$, it holds
$$
g_V^{\bar x}(xy) g_I^{\bar x y}(u) = g_V^{\bar x}(uy) g_I^{\bar x y}(x).
$$
\end{prop}
\begin{proof}
To prove the first statement, it is sufficient to recall that for any $I \subset V \in W$, $\Lambda \in W(V^c)$ and $\bar x \in X^{V^c}$, it holds
$$
g_I^{\bar x_\Lambda y}(x) = \frac{g_V^{\bar x_\Lambda}(xy)}{ \sum \limits_{u \in X^V} g_V^{\bar x_\Lambda}(uy)}, \qquad x,u \in X^I, y \in X^{V \backslash I}.
$$
The second statement of the proposition follows after the passage to the limit with respect to filtration~$\cm{F} = \{\cm{B}^{\Lambda_n}\}_{n \ge 1}$ in the expression
$$
g_V^{\bar x_{\Lambda_n \backslash V}}(xy) g_I^{\bar x_{\Lambda_n \backslash I} y}(u) = g_V^{\bar x_{\Lambda_n \backslash V}}(uy) g_I^{\bar x_{\Lambda_n \backslash I} y}(x), \qquad x,u \in X^I, y \in X^{V \backslash I}, \bar x \in X^{V^c},
$$
which is obvious for finite-conditional distributions.
\end{proof}

We will say that the system $G(P,\cm{F})$ is (\emph{strictly}) \emph{positive} if for all $V \in W$ and $\bar x \in \cm{X}_V(\cm{F})$, one has $g_V^{\bar x}(x)>0$, $x \in X^V$. The system $G(P,\cm{F})$ will be called \emph{quasilocal} if its elements are quasilocal as functions on the set of admissible boundary conditions, that is, if for any $V \in W$ and $x \in X^V$, it holds
$$
\lim \limits_{\Lambda \uparrow V^c} \sup
\limits_{{ \bar x, \bar y \in \cm{X}_V(\cm{F}): \bar x_\Lambda = \bar y_\Lambda}}
\left| g_V^{\bar x}(x) - g_V^{\bar y}(x) \right| = 0.
$$

The system $G(P,\cm{F})$ does not restore the random field $P$ since there may exists a random field $\hat P$, different from $P$, such that $G(\hat P,\cm{F}) = G(P,\cm{F})$ (the corresponding example will be given in the next subsection). In this regard, the question arises under what conditions the system $G(P,\cm{F})$ uniquely determines random field $P$. Using arguments similar to the proof of Theorem 4 in~\cite{Dobr}, one can show that if a random field $P$ satisfies the uniform strong mixing condition and has a positive and quasilocal system of conditional distributions $G(P,\cm{F})$ (corresponding to some filtration $\cm{F}$), then $G(P, \cm{F})$ uniquely determines $P$.

\begin{remark}
\label{KolmCond}
It is customary to use Kolmogorov's descriptive definition of conditional distribution of random field $P = \{P_V, V \in W\}$, which forms a set $Q(P) = \{q_V^{\bar x}, \bar x \in X^{V^c}, V \in W\}$ of probability distributions $q_V^{\bar x}$ on $X^V$ defined for all possible boundary conditions $\bar x \in X^{V^c}$.

It can be shown that if $G(P,\cm{F}) = \{g_V^{\bar x}, \bar x \in \cm{X}_V(\cm{F}), V \in W\}$ is a system of conditional distributions of the random field $P$ with respect to some filtration $\cm{F}$, then there exists a version $Q(P) = \{q_V^{\bar x}, \bar x \in X^{V^c}, V \in W\}$ of conditional distribution of~$P$ such that for all $V \in W$ and $\bar x \in \cm{X}_V(\cm{F})$, it holds
$$
q_V^{\bar x}(x) = g_V^{\bar x}(x), \qquad x \in X^V.
$$
$\bigtriangleup$
\end{remark}

\subsection{One-point conditional distribution}

A special role in the further presentation will be played by the systems $G_1(P,\cm{F}) = \{g_t^{\bar x}, \bar x \in \cm{X}_t(\cm{F}), t \in \mathbb{Z}^d\}$ of \emph{one-point conditional distributions} of the random field $P$ with respect to filtration $\cm{F}$. The connection between the elements of $G_1(P,\cm{F})$ is presented in the following statement, the proof of which is similar to the proof of Proposition~\ref{admiss}.

\begin{prop}
\label{admiss-1}
Let $G_1(P, \cm{F})$ be a system of one-point conditional distributions of the random field $P$ with respect to filtration $\cm{F}$. Then for any $t,s \in \mathbb{Z}^d$, $x,u \in X^t$, $y,v \in X^s$ and $\bar x \in \cm{X}_{\{t,s\}}(\cm{F})$, it holds
$$
g_t^{\bar x y}(x) g_s^{\bar x x}(v) g_t^{\bar x v}(u) g_s^{\bar x u}(y) = g_t^{\bar x y}(u) g_s^{\bar x u}(v) g_t^{\bar x v}(x) g_s^{\bar x x}(y).
$$
\end{prop}

The following statement holds true (see also~\cite{Kh22} and~\cite{DN04}).

\begin{T}
\label{GP1->GP}
For any random field $P$ and filtration $\cm{F}$, positive system $G_1(P,\cm{F})$ restores $G(P,\cm{F})$.
\end{T}
\begin{proof}
Consider some finite set $V \in W$ and let $\bar x \in \cm{X}_V(\cm{F})$ be an admissible boundary condition with respect to $\cm{F}$. By Proposition~\ref{admiss} for any $t \in V$ and $y \in X^{V \backslash t}$, one-point conditional distributions $g_t^{\bar x y}$ are well-defined. It remains to note that the following explicit formula expressing the elements of $G(P,\cm{F})$ in terms of the elements of $G_1(P,\cm{F})$ is valid:
\begin{equation}
\label{Q1P->QP}
g_V^{\bar x}(x) = \prod \limits_{j=1}^n \frac{g_{t_j}^{\bar x (xu)_j} (x_{t_j})}{g_{t_j}^{\bar x (xu)_j} (u_{t_j})} \cdot \left( \sum \limits_{z \in X^V} \prod \limits_{j=1}^n \frac{g_{t_j}^{\bar x (z u)_j} (z_{t_j})}{g_{t_j}^{\bar x (z u)_j} (u_{t_j})} \right)^{-1}, \quad x \in X^V, \bar x \in \cm{X}_V(\cm{F}),
\end{equation}
where $u \in X^V$ is an arbitrary configuration and
$$
\begin{array}{l}
(xu)_j = x_{t_j} ... x_{t_{j-1}} u_{t_{j+1}} ... u_{t_n}, \quad 1 < j < n, \qquad
(xu)_1 = u_{t_2} ... u_{t_n}, \quad (xu)_n = x_{t_1} x_{t_2} ... x_{t_{n-1}}.
\end{array}
$$
\end{proof}

By virtue of~\eqref{Q1P->QP}, many properties of conditional distributions $G(P,\cm{F})$ are consequences of the properties of the elements of its subsystem $G_1(P,\cm{F})$. For example, these are the properties of positivity, homogeneity and quasilocality.

An important class of random fields is formed by Markov random fields, the definition of which is formulated in terms of one-point conditional probabilities. A random field $P$ is called a \emph{Markov random field} (with respect to the neighborhood system $\partial$) if its one-point conditional distributions $G_1(P,\cm{F})$ posses the following (Markov) property: for all $t \in \mathbb{Z} ^d$ and $\bar x \in \cm{X}_t(\cm{F})$, it holds
$$
g_t^{\bar x}(x) = g_t^{\bar x_{\partial t}} (x), \qquad x \in X^t.
$$
The equivalent condition for a random field to be Markovian is as follows: for any $t \in \mathbb{Z}^d$ and $\Lambda \in W(t^c)$ such that $\partial t \subset \Lambda$, if $z,y \in X^\Lambda$ coincide on $\partial t$, that is, $z_{\partial t} = y_{\partial t}$, then
$$
g_t^z(x) =  g_t^y (x), \qquad x,u \in X^t.
$$
Note that for a Markov random field $P$, sets of admissible boundary conditions with respect to any filtration $\cm{F}$ are $\cm{X}_t(\cm{F}) = X^{t^c}$, $t \in \mathbb{Z}^d$, and thus, all systems of conditional distributions constructed with respect to different filtrations coincide.

The system $G_1(P,\cm{F})$ of one-point conditional distributions does not restore the random field $P$. Let us bring the corresponding example (see also~\cite{DN19} and~\cite{Georgii}).

\noindent \textit{Example 1}. Let $T = \{1,2,3,...\}$ and let numbers $c_j$, $0 < c_j < 1$, $j \in T$, be such that $\prod \limits_{j \in T} c_j >0$. Put $k_t = \prod \limits_{j=t}^\infty c_j$, $t \in T$. Let $P^{+}$ and $P^{-}$ be Markov random fields on $T$ with the same state space $X = \{ -1,1 \}$, initial distribution and transition probabilities of which are defined as
$$
p_1^{+} (x_1) = \frac{1 + x_1 k_1}{2}, \qquad p_1^{-} (x_1) = \frac{1 - x_1 k_1}{2},
$$
and
$$
p_t^{+} (x_{t-1}, x_t) = \frac{1 + c_{t-1} x_{t-1} x_t}{2} \cdot \frac{1 + x_t k_t}{1 + x_{t-1} k_{t-1}}, \quad p_t^{-} (x_{t-1}, x_t) = \frac{1 + c_{t-1} x_{t-1} x_t}{2} \cdot \frac{1 - x_t k_t}{1 - x_{t-1} k_{t-1}},
$$
respectively, where $x_t \in X^t$, $t \in T$. Let us show that the random fields $P^{+}$ and $P^{-}$ have the same system of one-point conditional distributions (with respect to any filtration $\cm{F}$).

It is not difficult to see that for all $n \in T$, finite-dimensional distributions $P_{\{1,2,...,n\}}^+$ and $P_{\{1,2,...,n\}}^-$ of the random processes $P^{+}$ and $P^{-}$ have the following form:
$$
P_{\{1,2,...,n\}}^\pm (x_1 x_2 ... x_n) = \left( \prod \limits_{j=1}^{n - 1} \frac{1 + c_j x_j x_{j+1} }{2} \right) \cdot \frac{1 \pm x_n k_n}{2}, \qquad x_j \in X^j, 1 \le j \le n.
$$
Let ${}^+ q_t^y$ and ${}^- q_t^y$, $y \in X^\Lambda$, $\Lambda \in W(t^c)$, $t \in T$, be one-point finite-conditional distributions of the random fields $P^{+}$ and $P^{-}$, respectively. It can be shown by direct computations that for all $t > 1$ and $\Lambda = \{1,2,...,t-1,t+1,...,t+n\}$, $n \in \mathbb{N}$, it holds
$$
{}^+ q_t^y (x_t) = {}^- q_t^y (x_t) = \dfrac{(1 + c_{t-1} y_{t-1} x_t )(1+ c_t x_t y_{t+1}) }{2 (1 + c_{t-1} c_t y_{t-1} y_{t+1})}, \qquad x \in X^t, y \in X^\Lambda,
$$
and for $t = 1$,
$$
{}^+ g_1^y (x_1) = {}^- g_1^y (x_1) = \dfrac{1 + c_1 x_1 y_2}{2}, \quad x_1 \in X^1, \, y \in X^{\{2,3,...,n\}}.
$$
From here it follows that for any increasing sequence $\{\Lambda_n\}_{n \ge 1}$ of finite sets converging to $T$ and for all $\bar y \in X^{T \backslash t}$, the limits
$$
g_t^{\bar y}(x_t) = \lim \limits_{\Lambda_n \uparrow T} {}^\pm g_t^{\bar y_{\Lambda_n \backslash t}} (x_t) = \dfrac{(1 + c_{t - 1} \bar y_{t -1} x_t )(1+ c_t x_t \bar y_{t+1}) }{2 (1 + c_{t-1} c_t \bar y_{t-1} \bar y_{t+1})}, \qquad x_t \in X^{t}, \, t > 1,
$$
$$
g_1^{\bar y}(x_1) = \lim \limits_{\Lambda_n \uparrow T} {}^\pm g_1^{\bar y_{\Lambda_n \backslash \{1\}}} (x_1) = \dfrac{1 + c_1 x_1 \bar y_2}{2}, \qquad x_1 \in X^{\{1\}}.
$$
exist. Therefore, for any filtration $\cm{F}$, we have $G_1(P^{+},\cm{F}) = \{g_t^{\bar y}, \bar y \in X^{T \backslash t}, t \in T\} =G_1(P^{-},\cm{F})$. \qed \\

\subsection{Transition energy and Hamiltonian of a random field}

In this section, we introduce concepts of transition energy and the Hamiltonian of a random field which are important for further presentation and, possibly, for the theory of random fields itself. Their autonomous (without appeal to the notion of a random field) analogues were introduced in~\cite{DN19}.

Let the system of conditional distributions $G(P,\cm{F}) = \{g_V^{\bar x}, \bar x \in \cm{X}_V(\cm{F}), V \in W\}$ of the random field $P$ with respect to some filtration $\cm{F}$ be positive. Denote by $\boldsymbol{\Delta}(P,\cm{F}) = \{ \Delta_V^{\bar x}, \bar x \in \cm{X}_V(\cm{F}), V \in W \}$ the set of functions defined by
\begin{equation}
\label{DeltaP}
\Delta_V^{\bar x}(x,u) = \ln \frac{g_V^{\bar x}(x)}{g_V^{\bar x}(u)}, \qquad x,u \in X^V, \bar x \in \cm{X}_V(\cm{F}), V \in W.
\end{equation}
It is not difficult to see that for all $V \in W$ and $\bar x \in \cm{X}_V(\cm{F})$, the following equality takes place:
\begin{equation}
\label{DeltaP-1}
\Delta_V^{\bar x}(x,u) = \Delta_V^{\bar x}(x,y) + \Delta_V^{\bar x}(y,u), \qquad x,u,y \in X^V.
\end{equation}
In particular, from here it follows that
$$
\Delta_V^{\bar x}(x,u) = - \Delta_V^{\bar x}(u,x), \qquad \Delta_V^{\bar x}(x,x) = 0, \qquad x,u \in X^V.
$$
The set $\boldsymbol{\Delta}(P,\cm{F})$ will be called the \emph{system of transition energies of the random field $P$ with respect to filtration $\cm{F}$}.

Note that the elements of $\boldsymbol{\Delta}(P,\cm{F})$ can be computed directly by the finite-dimensional distribution of $P$ as
$$
\Delta_V^{\bar x}(x,u) = \lim \limits_{n \to \infty} \ln \frac{P_{\Lambda_n}(x \bar x_{\Lambda_n \backslash V})}{P_{\Lambda_n}(u \bar x_{\Lambda_n  \backslash V})}, \qquad x,u \in X^V, \bar x \in \cm{X}_V(\cm{F}), V \in W,
$$
where $\cm{F} = \{\cm{B}^{\Lambda_n}\}_{n \ge 1}$.

The next statement follows from the definition of $\boldsymbol{\Delta}(P,\cm{F})$.

\begin{T}
\label{prop-QP-Gibbs}
Let $P$ be a random field having a positive system $G(P,\cm{F}) = \{ g_V^{\bar x}, \bar x \in \cm{X}_V(\cm{F}), V \in W \}$ of conditional distributions with respect to some filtration $\cm{F}$. Then the elements of $G(P,\cm{F})$ necessarily have a Gibbs form
$$
g_V^{\bar x} (x) = \frac{ \exp \{\Delta_V^{\bar x} (x,u)\} }{\sum \limits_{\alpha \in X^V} \exp \{\Delta_V^{\bar x} (\alpha,u)\}}, \qquad x \in X^V, \bar x \in \cm{X}_V(\cm{F}), V \in W,
$$
where $u \in X^V$ and $\boldsymbol{\Delta}(P,\cm{F}) = \{ \Delta_V^{\bar x}, \bar x \in \cm{X}_V(\cm{F}), V \in W \}$ is the system of transition energies of $P$ with respect to $\cm{F}$.
\end{T}

Let us also mention the following property of the system $\boldsymbol{\Delta}(P,\cm{F})$: for all disjoint $V, I \in W$ and $\bar x \in \cm{X}_{V \cup I}(\cm{F})$, it holds
\begin{equation}
\label{DeltaP-2}
\Delta_{V \cup I}^{\bar x}(xy,uv) = \Delta_V^{\bar x y}(x,u) + \Delta_I^{\bar x x}(y,v), \qquad x,u \in X^V, y,v \in X^I.
\end{equation}
Indeed, by Proposition~\ref{admiss}, boundary conditions $\bar x y$ and $\bar x x$ are also admissible. It remains to note that for all pairwise disjoint sets $V, I, \Lambda \in W$ and configurations $x,u \in X^V$, $y,v \in X^I$, $z \in X^\Lambda$, the following relation is satisfied:
$$
\ln \frac{P_{V \cup I \cup \Lambda}(x y z)}{P_{V \cup I \cup \Lambda}(u v z)} = \ln \frac{P_{V \cup I \cup \Lambda}(x y z)}{P_{V \cup I \cup \Lambda}(u y z)} + \ln \frac{P_{V \cup I \cup \Lambda}(u y z)}{P_{V \cup I \cup \Lambda}(u v z)}.
$$

Since the elements of $\boldsymbol{\Delta}(P,\cm{F})$ satisfy~\eqref{DeltaP-1}, for all $V \in W$ and $\bar x \in \cm{X}_V(\cm{F})$, they can be presented in the form
\begin{equation}
\label{Delta->HP}
\Delta_V^{\bar x}(x,u) = H_V^{\bar x}(u) - H_V^{\bar x}(x), \qquad x,u \in X^V,
\end{equation}
where function $H_V^{\bar x}$ can be interpreted as a Hamiltonian of the random field $P$ in finite volume $V$ under admissible boundary condition $\bar x$. The set $H(P,\cm{F}) = \{H_V^{\bar x}, \bar x \in \cm{X}_V(\cm{F}), V \in W\}$ of such functions will be called a \emph{Hamiltonian corresponding to the random field $P$ with respect to filtration $\cm{F}$}.

It is clear, that for each $V \in W$ and $\bar x \in \cm{X}_V(\cm{F})$, there exist different functions $H_V^{\bar x}$ on $X^V$, satisfying~\eqref{Delta->HP}. In particular, according to~\eqref{DeltaP}, one can put
$$
H_V^{\bar x}(x) = - \ln q_V^{\bar x}(x) + c_V(\bar x), \qquad x \in X^V,
$$
where $c_V$ is some function on $X^{V^c}$, $V \in W$. Thus, to any random field $P$ and any filtration $\cm{F}$ it may correspond various Hamiltonians $H(P,\cm{F})$.

A direct corollary of Theorem~\ref{prop-QP-Gibbs} is the following result.

\begin{corol}
\label{corol-H}
Let the system $G(P,\cm{F})$ of conditional distributions of the random field $P$ with respect to some filtration $\cm{F}$ be positive. Then its elements have a Gibbs form
$$
g_V^{\bar x} (x) = \frac{ \exp \{- H_V^{\bar x} (x)\} }{\sum \limits_{\alpha \in X^V} {\exp \{ -H_V^{\bar x} (\alpha)\} }}, \qquad x \in X^V, \bar x \in \cm{X}_V(\cm{F}), V \in W,
$$
where $H(P,\cm{F}) = \{H_V^{\bar x}, \bar x \in \cm{X}_V(\cm{F}), V \in W\}$ is a Hamiltonian corresponding to the random field $P$ with respect to filtration $\cm{F}$.
\end{corol}

Let us also note that from~\eqref{DeltaP-2} it follows that the elements of the Hamiltonian $H(P,\cm{F}) = \{H_V^{\bar x}, \bar x \in \cm{X}_V(\cm{F}), V \in W\}$ are consistent in the following sense: for all disjoint $V, I \in W$ and $\bar x \in \cm{X}_{V \cup I}(\cm{F})$, it holds
$$
H_{V \cup I}^{\bar x}(xy) + H_V^{\bar x y}(u) = H_{V \cup I}^{\bar x}(uy) + H_V^{\bar x y}(x), \qquad x,u \in X^V, y,v \in X^I.
$$

Now, let us consider the set $\boldsymbol{\Delta}_1(P,\cm{F}) = \{ \Delta_t^{\bar x}, \bar x \in \cm{X}_t(\cm{F}), t \in \mathbb{Z}^d \}$ representing the \emph{system of one-point transition energies of the random field $P$ with respect to filtration $\cm{F}$}. It is clear that for all $t \in \mathbb{Z}^d$ and $\bar x \in \cm{X}_t(\cm{F})$, one-point transition energies satisfy the relation
$$
\Delta_t^{\bar x}(x,u) = \Delta_t^{\bar x}(x,y) + \Delta_t^{\bar x}(y,u), \qquad x,u,y \in X^t.
$$
Moreover, it is not difficult to verify that the elements of $\boldsymbol{\Delta}_1(P,\cm{F})$ are consistent in the following sense: for all $t,s \in \mathbb{Z}^d$ and $\bar x \in \cm{X}_{\{t,s\}}(\cm{F})$, it holds
\begin{equation}
\label{Delta1P-2}
\Delta_t^{\bar x y}(x,u) + \Delta_s^{\bar x u}(y,v) = \Delta_s^{\bar x x}(y,v) + \Delta_t^{\bar x v}(x,u), \qquad x,u \in X^t, y,v \in X^s.
\end{equation}

The system $\boldsymbol{\Delta}_1(P,\cm{F})$ of one-point transition energies will be called \emph{quasilocal} (with respect to boundary conditions) if for all $t \in \mathbb{Z}^d$ and $x,u \in X^t$, one has
$$
\lim \limits_{\Lambda \uparrow t^c} \sup \limits_{ \bar x, \bar y \in \cm{X}_t(\cm{F}): \bar x_\Lambda = \bar y_\Lambda} \left| \Delta_t^{\bar x}(x,u) - \Delta_t^{\bar y}(x,u) \right| = 0.
$$

As in the case of $\boldsymbol{\Delta}(P,\cm{F})$, the elements of $\boldsymbol{\Delta}_1(P,\cm{F})$ can be represented in the form
$$
\Delta_t^{\bar x}(x,u) = H_t^{\bar x}(u) - H_t^{\bar x}(x), \qquad x,u \in X^t, \bar x \in \cm{X}_t(\cm{F}), t \in \mathbb{Z}^d,
$$
and the set of functions $H_1(P,\cm{F}) = \{H_t^{\bar x}, \bar x \in \cm{X}_t(\cm{F}), t \in \mathbb{Z}^d\}$ can be interpreted as a \emph{one-point Hamiltonian corresponding to the random field $P$ with respect to filtration $\cm{F}$}. Due to~\eqref{Delta1P-2}, the elements of $H_1(P,\cm{F})$ satisfy the following consistency conditions: for all $t,s \in \mathbb{Z}^d$ and $\bar x \in \cm{X}_{\{t,s\}}(\cm{F})$, $x,u \in X^t$, $y,v \in X^s$, it holds
$$
H_t^{\bar x y}(x) + H_s^{\bar x x}(v) + H_t^{\bar x v}(u) + H_s^{\bar x u}(y) = H_t^{\bar x y}(u) + H_s^{\bar x u}(v) + H_t^{\bar x v}(x) + H_s^{\bar x x}(y).
$$

Note that, unlike the transition energies, the Hamiltonian does not have the property of additivity.

\subsection{Specification and transition energy field}
\label{Sect-auto}

Along with the considered systems of functions (systems of conditional distributions and transition energies, as well as Hamiltonian), constructed from a random field, we will use their autonomously (\textit{a priori}) given analogues.

A set $Q=\{q_V^{\bar x}, \bar x \in X^{V^c}, V \in W\}$ of probability distributions $q_V^{\bar x}$ on $X^V$ parameterized by infinite boundary conditions $\bar x$ is called a \emph{specification} if for all $I \subset V \in W$ and $\bar x \in X^{V^c}$, its elements satisfy the following consistency conditions:
$$
q_V^{\bar x}(xy) q_I^{\bar x y}(u) = q_V^{\bar x}(uy) q_I^{\bar x y}(x), \qquad x,u \in X^I, y \in X^{V \backslash I}.
$$

A set $Q_1=\{q_t^{\bar x}, \bar x \in X^{t^c}, t \in \mathbb{Z}^d\}$ of positive one-point probability distributions $q_t^{\bar x}$ on $X^t$ parameterized by infinite boundary conditions $\bar x$ is called a \emph{1--specification} (see~\cite{DN04}) if for all $t,s \in \mathbb{Z}^d$ and $\bar x \in X^{{\{t,s\}}^c}$, it holds
$$
q_t^{\bar x y}(x) q_s^{\bar x x}(v) q_t^{\bar x v}(u) q_s^{\bar x u}(y) = q_t^{\bar x y}(u) q_s^{\bar x u}(v) q_t^{\bar x v}(x) q_s^{\bar x x}(y), \qquad x,u \in X^t, y,v \in X^s.
$$

In~\cite{DN04}, it was shown, in particular, that any 1--specification $Q_1$ uniquely determines the specification $Q$ of which it is a subsystem.

Autonomously defined analogues of the system of transition energies and the system of one-point transition energies were introduced in~\cite{DN19}.

A set $\boldsymbol{\Delta} = \{ \delta_V^{\bar x}, \bar x \in X^{V^c}, V \in W \}$ of functions $\delta_V^{\bar x}(x,u)$, $x,u \in X^V$, is called \emph{transition energy field} if its elements satisfy the following consistency conditions: for all $V \in W$ and $\bar x \in X^{V^c}$, it holds
$$
\delta_V^{\bar x}(x,u) = \delta_V^{\bar x}(x,y) + \delta_V^{\bar x}(y,u), \qquad x,u,y \in X^V;
$$
and for all disjoint $V, I \in W$ and $\bar x \in X^{(V \cup I)^c}$,
$$
\delta_{V \cup I}^{\bar x}(xy,uv) = \delta_V^{\bar x y}(x,u) + \delta_I^{\bar x x}(y,v), \qquad x,u \in X^V, y,v \in X^I.
$$

A set $\boldsymbol{\Delta}_1 = \{ \delta_t^{\bar x}, \bar x \in X^{t^c}, t \in \mathbb{Z}^d \}$ of functions $\delta_t^{\bar x}(x,u)$, $x,u \in X^t$, is called \emph{one-point transition energy field} if its elements satisfy the following consistency conditions: for all $t \in \mathbb{Z}^d$ and $\bar x \in X^{t^c}$, it holds
$$
\delta_t^{\bar x}(x,u) = \delta_t^{\bar x}(x,y) + \delta_t^{\bar x}(y,u), \qquad x,u,y \in X^t;
$$
and for all $t,s \in \mathbb{Z}^d$ and $\bar x \in X^{\{t,s\}^c}$,
$$
\delta_t^{\bar x y}(x,u) + \delta_s^{\bar x u}(y,v) = \delta_s^{\bar x x}(y,v) + \delta_t^{\bar x v}(x,u), \qquad x,u \in X^t, y,v \in X^s.
$$

In~\cite{DN19}, a fundamental connection was established between the 1-specification and the one-point transition energy field according to which a system $Q_1=\{q_t^{\bar x}, \bar x \in X^{t^c}, t \in \mathbb{Z}^d\}$ is a 1--specification if and only if its elements are representable in a Gibbs form
$$
q_t^{\bar x}(x) = \frac{ \exp \{\delta_t^{\bar x} (x,u)\} }{\sum \limits_{z \in X^t} \exp \{\delta_t^{\bar x} (\alpha,u)\} }, \qquad x \in X^t, \bar x \in X^{t^c}, t \in \mathbb{Z}^d,
$$
where $\boldsymbol{\Delta}_1 = \{ \delta_t^{\bar x}, \bar x \in X^{t^c}, t \in \mathbb{Z}^d \}$ is a one-point transition energy field. A similar result holds for the specification and the transition energy field.

Finally, let us note that in~\cite{DN19}, the axiomatic definition of a (one-point) Hamiltonian was proposed.

A \emph{Hamiltonian} is a set $H = \{H_V^{\bar x}, \bar x \in X^{V^c}, V \in W\}$ of functions on $X^V$, $V \in W$, consistent in the following sense: for all disjoint $V, I \in W$ and $\bar x \in X^{(V \cup I)^c}$, it holds
$$
H_{V \cup I}^{\bar x}(xy) + H_V^{\bar x y}(u) = H_{V \cup I}^{\bar x}(uy) + H_V^{\bar x y}(x), \qquad x,u \in X^V, y,v \in X^I.
$$

A \emph{one-point Hamiltonian} is a set $H_1 = \left\{H_t^{\bar x}, \bar x \in X^{t^c}, t \in \mathbb{Z}^d \right\}$ of functions on $X^t$, $t \in \mathbb{Z}^d$, satisfying the following consistency conditions: for all $t,s \in \mathbb{Z}^d$ and $\bar x \in X^{\{t,s\}^c}$, it holds
$$
H_t^{\bar x y}(x) + H_s^{\bar x x}(v) + H_t^{\bar x v}(u) + H_s^{\bar x u}(y) = H_t^{\bar x y}(u) + H_s^{\bar x u}(v) + H_t^{\bar x v}(x) + H_s^{\bar x x}(y), \qquad x,u \in X^t, y,v \in X^s.
$$

It is not difficult to see that if $H_1 = \left\{H_t^{\bar x}, \bar x \in X^{t^c}, t \in \mathbb{Z}^d \right\}$ is a one-point Hamiltonian, then the set $\boldsymbol{\Delta}_1 = \{ \delta_t^{\bar x}, \bar x \in X^{t^c}, t \in \mathbb{Z}^d \}$ of functions
\begin{equation}
\label{Delta-H}
\delta_t^{\bar x}(x,u) = H_t^{\bar x}(u) - H_t^{\bar x}(x), \qquad x,u \in X^t,
\end{equation}
forms a one-point transition energy field. Conversely, the elements of any one-point transition energy field $\boldsymbol{\Delta}_1$ can be represented in the form~\eqref{Delta-H} in terms of some one-point Hamiltonian $H_1$. Moreover, a system $Q_1$ is a 1--specification if and only if its elements are representable in a Gibbs form
$$
q_t^{\bar x}(x) = \frac{ \exp \{ -H_t^{\bar x} (x)\} }{\sum \limits_{\alpha \in X^t} \exp \{-H_t^{\bar x} (\alpha)\} }, \qquad x \in X^t, \bar x \in X^{t^c}, t \in \mathbb{Z}^d,
$$
with some one-point Hamiltonian $H_1$.

The basic example of Hamiltonian is a Hamiltonian $H^\Phi$ constructed by convergent potential $\Phi$. A specification $Q^\Phi$ corresponding to $H^\Phi$ is called a (Gibbs) specification constructed by $\Phi$.

\begin{remark}
In the DLR--framework, to show that a specification is Gibbsian, it is necessary to find the corresponding potential. For this, the Mobius formula is usually used (this approach was first proposed by Grimmett~\cite{Grimmet}), which leads, as a rule, to a very complicated expression for the potential. Moreover, as Gandolfi and Lenarda~\cite{GL} note, when using the Mobius formula, the explicit dependence of the potential from the spin values of physical systems is lost.


The results of~\cite{DN19} show that for a given specification there is no need to look for a potential since any specification admits a Gibbs representation in terms of the transition energy field. The latter is uniquely determined by the specification, and, in particular, the translational invariance of the specification is equivalent to the translational invariance of the corresponding transition energy field. At the same time, there is no such uniqueness for the potential corresponding to the specification. In addition, as shown in~\cite{BGMMT}, there are examples of translation-invariant specifications, representable in Gibbs form with a uniformly convergent potential, for which there is no translation-invariant potential. $\quad \bigtriangleup$
\end{remark}

\section{Probabilistic definition of a Gibbs random field and outlines of the corresponding theory}
\label{sect-Gibbs}

In this section, we give a probabilistic definition of a Gibbs random field and show how the structure of the theory of Gibbs random fields, built using this definition, should look like. The circle of problems that are natural to consider within the framework of this theory is described. Also, we give a comparative analysis of the proposed definition with the DLR--definition and other existing definitions of a Gibbs random field.

\subsection{Probabilistic definition}

As it was mentioned above, by the probabilistic definition of a Gibbs random field we mean a definition that is based solely on its intrinsic properties, namely, on the properties of its finite-dimensional distributions. Such a definition should be free from any physical concepts, for example, from the concept of potential.

There are definitions of a Gibbs random field that do not use the concept of potential but, nevertheless, are not probabilistic since all such definitions are always based on an \textit{a priori} given object --- a specification (see, for example,~\cite{Fernandez, LeNy, MaesKo}).

The first purely probabilistic definition of a Gibbs random field was given in~\cite{DN09}. Below, we present a probabilistic definition of a Gibbs random field, which is more suitable for constructing the theory of Gibbs random fields on its basis.

\noindent \textit{Definition}.
A random field $P = \{P_V, V \in W\}$ will be called a \emph{Gibbs random field} if there exists a filtration $\cm{F} = \{\cm{B}^{\Lambda_n}\}_{n \ge 1}$ and sets of admissible boundary conditions $\cm{Y}_t(\cm{F}) \subset \cm{X}_t(\cm{F})$, $P(\cm{Y}_t(\cm{F})) = 1$, $t \in \mathbb{Z}^d$, such that for all $t \in \mathbb{Z}^d$ and any $\bar x \in \cm{Y}_t(\cm{F})$, the limits of finite-conditional distributions
\begin{equation}
\label{q1-Gibbs}
\lim \limits_{n \to \infty} \frac{P_{\Lambda_n}(x \bar x_{\Lambda_n \backslash t})}{P_{\Lambda_n \backslash t} (\bar x_{\Lambda_n \backslash t})} = \lim \limits_{n \to \infty} g_t^{\bar x_{\Lambda_n \backslash t}}(x) = g_t^{\bar x} (x), \qquad x \in X^t,
\end{equation}
are positive and convergence with respect to $\bar x \in \cm{Y}_t(\cm{F})$ is uniform.

We will call this definition of Gibbs random fields the $\Pr$--definition. Filtration $\cm{F}$ will be called a \emph{determining filtration}. It is clear that any Gibbs random field may have several determining filtrations.

The following theorem states the behavior of finite-conditional distributions of a Gibbs random field with respect to any filtration.

\begin{T}
\label{GibbsProperty}
Let $P$ be a Gibbs random field, $\cm{F} = \{\cm{B}^{\Lambda_n}\}_{n \ge 1}$ be its determining  filtration and $\cm{Y}_t(\cm{F})$, $t \in \mathbb{Z}^d$, be corresponding sets of admissible boundary conditions. Then for any filtration $\cm{F}' = \{\cm{B}^{\Lambda'_n}\}_{n \ge 1}$ and for any $t \in \mathbb{Z}^d$, $x \in X^t$, the sequence of finite-conditional distributions $\left\{g_t^{\bar x_{\Lambda'_n \backslash t}}(x)\right\}_{n \ge 1}$ converge (when $n \to \infty$) to the same limit~\eqref{q1-Gibbs} uniformly with respect to $\bar x \in \cm{Y}_t(\cm{F})$.
\end{T}
\begin{proof}
First let us show that one-point conditional distributions of Gibbs random field $P$ with respect to its basic filtration $\cm{F}$ are quasilocal on the corresponding sets of admissible boundary conditions $\cm{Y}_t(\cm{F})$, $P(\cm{Y}_t(\cm{F})) = 1$, $t \in \mathbb{Z}^d$, that is,
\begin{equation}
\label{quasilocY}
\lim \limits_{\Lambda \uparrow t^c} \sup
\limits_{{ \bar x, \bar y \in \cm{Y}_t(\cm{F}): \bar x_\Lambda = \bar y_\Lambda}}
\left| g_t^{\bar x}(x) - g_t^{\bar y}(x) \right| = 0, \qquad x \in X^t.
\end{equation}
Indeed, according to $\Pr$--definition, for any $\varepsilon > 0$ and sufficiently large $n$, for any $x \in X^t$, it holds
$$
\sup \limits_{\bar x \in \cm{Y}_t(\cm{F})} \left| g_t^{\bar x}(x) - g_t^{\bar x_{\Lambda_n \backslash t}}(x) \right| \le \varepsilon.
$$
Therefore, for any $\Lambda \in W(t^c)$ such that $\Lambda_n \backslash t \subset \Lambda$, we obtain
$$
\begin{array}{l}
\sup \limits_{\bar x, \bar y \in \cm{Y}_t(\cm{F}): \bar x_\Lambda = \bar y_\Lambda} \vert g_t^{\bar x}(x) - g_t^{\bar y}(x) \vert \le \sup \limits_{\bar x \in \cm{Y}_t(\cm{F})} \left| g_t^{\bar x}(x) - g_t^{\bar x_{\Lambda_n \backslash t}}(x) \right| + \sup \limits_{\bar y \in \cm{Y}_t(\cm{F})} \left| g_t^{\bar y}(x) - g_t^{\bar y_{\Lambda_n \backslash t}}(x) \right|
\le 2\varepsilon.
\end{array}
$$

Further, for all $t \in \mathbb{Z}^d$, $\Lambda \in W(t^c)$ and $x \in X^t$, $z \in X^\Lambda$, we have
$$
g_t^z(x) = \frac{1}{P_\Lambda(z)} \int \limits_{\{\bar y \in \cm{Y}_t(\cm{F}): \bar y_\Lambda = z\}} g_t^{\bar y}(x) P_{t^c}(d \bar y),
$$
and hence, the following inequalities hold
$$
\inf \limits_{\bar y \in \cm{Y}_t(\cm{F}): \bar y_\Lambda = z} g_t^{\bar y}(x) \le g_t^z(x) \le \sup \limits_{\bar y \in \cm{Y}_t(\cm{F}): \bar y_\Lambda = z} g_t^{\bar y}(x).
$$
Thus, for any $\bar x \in \cm{Y}_t(\cm{F})$, we can write
$$
\begin{array}{l}
\left| g_t^{\bar x} (x) - g_t^{\bar x_\Lambda} (x) \right| \le \max \left\{ \left| g_t^{\bar x}(x) - \sup \limits_{\bar y \in \cm{Y}_t(\cm{F}): \bar y_\Lambda = \bar x_\Lambda} g_t^{\bar y}(x) \right|, \left| g_t^{\bar x}(x) - \inf \limits_{\bar y \in \cm{Y}_t(\cm{F}): \bar y_\Lambda = \bar x_\Lambda} g_t^{\bar y}(x) \right| \right\} \le \\
\\
\le \sup \limits_{\bar y \in \cm{Y}_t(\cm{F}): \bar y_\Lambda = \bar x_\Lambda} \left| g_t^{\bar x}(x) - g_t^{\bar y}(x) \right|.
\end{array}
$$

Now let $\cm{F}' = \{\cm{B}^{\Lambda'_n}\}_{n \ge 1}$ be some filtration. Then for any $n \ge 1$, we have
$$
\sup \limits_{\bar x \in \cm{Y}_t(\cm{F})} \left| g_t^{\bar x} (x) - g_t^{\bar x_{\Lambda'_n \backslash t}} (x) \right| \le \sup \limits_{\bar x \in \cm{Y}_t(\cm{F})} \sup \limits_{\bar y \in \cm{Y}_t(\cm{F}): \bar y_{\Lambda'_n \backslash t} = \bar x_{\Lambda'_n \backslash t}} \left| g_t^{\bar x}(x) - g_t^{\bar y}(x) \right|,
$$
and hence, by quasilocality~\eqref{quasilocY} of one-point conditional distributions on the set $\cm{Y}_t(\cm{F})$, we obtain
$$
\lim \limits_{\Lambda'_n \uparrow \mathbb{Z}^d} \sup \limits_{\bar x \in \cm{Y}_t(\cm{F})} \left| g_t^{\bar x}(x) - g_t^{\bar x_{\Lambda'_n \backslash t}}(x) \right| = 0
$$
for all $t \in \mathbb{Z}^d$ and $x \in X^t$.
\end{proof}

In the next proposition, we show the existence of universal sets of boundary conditions on which finite-conditional distributions of a Gibbs random field converge uniformly.

\begin{prop}
\label{prop-Yt}
Let $P$ be a Gibbs random field. Then for each $t \in \mathbb{Z}^d$, there exists a set of boundary conditions $\cm{Y}_t \subset X^{t^c}$ such that $P(\cm{Y}_t) = 1$ and the limits of finite-conditional distributions~\eqref{q1-Gibbs} exist for any filtration $\cm{F} = \{\cm{B}^{\Lambda_n}\}_{n \ge 1}$, are positive, and convergence with respect to $\bar x \in \cm{Y}_t$ is uniform.
\end{prop}
\begin{proof}
Let $P$ be a Gibbs random field, $\cm{F} = \{\cm{B}^{\Lambda_n}\}_{n \ge 1}$ be its determining  filtration and $\cm{Y}_t(\cm{F})$, $P(\cm{Y}_t(\cm{F}))=1$, $t \in \mathbb{Z}^d$, be corresponding sets of admissible boundary conditions. According to Theorem~\ref{GibbsProperty}, for any filtration $\cm{F}' = \{\cm{B}^{\Lambda'_n}\}_{n \ge 1}$ and for any $t \in \mathbb{Z}^d$, $x \in X^t$, the sequence of finite-conditional distributions $\left\{g_t^{\bar x_{\Lambda'_n \backslash t}}(x)\right\}_{n \ge 1}$ converge (when $n \to \infty$) uniformly with respect to $\bar x \in \cm{Y}_t(\cm{F})$ to positive limits. Hence, filtration $\cm{F}'$ is also determining for the Gibbs random field $P$. Since $\cm{Y}_t(\cm{F}) \subset \cm{X}_t(\cm{F'})$, $t \in \mathbb{Z}^d$, one can put $\cm{Y}_t(\cm{F'}) = \cm{Y}_t(\cm{F})$, $t \in \mathbb{Z}^d$.
\end{proof}

For a random field $P$ and any $t \in \mathbb{Z}^d$, consider $\cm{X}_t = \bigcap \limits_{\cm{F}} \cm{X}_t(\cm{F})$ where intersection is taken over all possible filtrations $\cm{F} = \{\cm{B}^{\Lambda_n}\}_{n \ge 1}$. Note that $\cm{X}_t \in \cm{B}^{t^c}$, and, in general, $P(\cm{X}_t)$ may differ from 1. However, for Gibbs random fields, we have the following result.

\begin{prop}
\label{prop-Yt-1}
If $P$ is a Gibbs random field, then for any $t \in \mathbb{Z}^d$, one has $P(\cm{X}_t) = 1$.
\end{prop}
\begin{proof}
According to Proposition~\ref{prop-Yt}, for the Gibbs random fields $P$, there exist sets of boundary conditions $\cm{Y}_t$, $t \in \mathbb{Z}^d$, on which its finite-conditional distributions converge with respect to any filtration $\cm{F}$. It remains to note that $\cm{Y}_t \subset \cm{X}_t$ and $P(\cm{Y}_t) = 1$.
\end{proof}

Proposition~\ref{prop-Yt-1} provides the following test for non-Gibbsianness of a random field.

\begin{prop}
\label{nonGibbs}
If for a random field $P$, $P(\cm{X}_t) \neq 1$ for some $t \in \mathbb{Z}^d$, then $P$ is not a Gibbs random field.
\end{prop}

Below, we present a criterion of Gibbsianness of random fields, which is formulated in terms of conditional distributions of the random field (cf. Criterion 3.4 in~\cite{DN09}).

\begin{T}
\label{crit-Gibbs-1}
A random field $P$ is a Gibbs random field if there exists a filtration $\cm{F}$ such that for each $t \in \mathbb{Z}^d$, there is a set of admissible boundary conditions $\cm{Y}_t(\cm{F}) \in \cm{X}_t(\cm{F})$, $P(\cm{Y}_t(\cm{F})) = 1$, on which one-point conditional distributions of $P$ with respect to $\cm{F}$ are quasilocal, and for each $\bar x \in \cm{Y}_t(\cm{F})$, one has $g_t^{\bar x}(x) > 0$, $x \in X^t$.

Conversely, if $P$ is a Gibbs random field, then there exist sets of boundary conditions $\cm{Y}_t \subset X^{t^c}$, $P(\cm{Y}_t) = 1$, $t \in \mathbb{Z}^d$, such that for any filtration $\cm{F}$ one-point conditional distributions $G_1(P,\cm{F})$ are quasilocal on $\cm{Y}_t$, $t \in \mathbb{Z}^d$, and for each $\bar x \in \cm{Y}_t$, one has $g_t^{\bar x}(x) > 0$, $x \in X^t$.
\end{T}

The proof of this theorem can be deduced from the proof of Theorem~\ref{GibbsProperty} and will therefore be omitted.

The following proposition establishes the Gibbsianness of a random filed in terms of its transition energies.

\begin{prop}
\label{prop-crit-Gibbs-2}
A random field $P$ is a Gibbs random field if there exists filtration $\cm{F}$ such that for all $t \in \mathbb{Z}^d$, there is a set of admissible boundary conditions $\cm{Y}_t(\cm{F})$, $P(\cm{Y}_t(\cm{F})) = 1$, on which transition energies of $P$ with respect to $\cm{F}$ are quasilocal, that is,
$$
\lim \limits_{\Lambda \uparrow t^c} \sup \limits_{ \bar x, \bar y \in \cm{Y}_t(\cm{F}): \bar x_\Lambda = \bar y_\Lambda} \left| \Delta_t^{\bar x}(x,u) - \Delta_t^{\bar y}(x,u) \right| = 0, \qquad x,u \in X^t.
$$
\end{prop}
\begin{proof}
First note that for all $t \in \mathbb{Z}^d$, one-point conditional distributions
$$
g_t^{\bar x}(x) = \frac{ \exp \{\Delta_t^{\bar x} (x,u)\} }{\sum \limits_{\alpha \in X^V} \exp \{\Delta_t^{\bar x} (\alpha,u)\}}, \qquad x \in X^t, \bar x \in \cm{Y}_t(\cm{F}),
$$
of $P$ with respect to $\cm{F}$ are positive and quasilocal on $\cm{Y}_t(\cm{F})$. It remains to apply Theorem~\ref{crit-Gibbs-1}.
\end{proof}

The converse statement that transition energies of a Gibbs random field are quasilocal can be proven under stricter conditions, namely, its conditional distributions have to be uniformly nonnull.

Next, we present properties of conditional distributions of Gibbs random fields, which directly follows from relation~\eqref{Q1P->QP}.

\begin{prop}
Let $P=\{P_V, V\in W\}$ be a Gibbs random field, $\cm{F} = \{\cm{B}^{\Lambda_n}\}_{n \ge 1}$ be its basic filtration and $\cm{Y}_t(\cm{F})$, $t \in \mathbb{Z}^d$, be corresponding sets of admissible boundary conditions. For all $V \in W$, put $\cm{Y}_V(\cm{F}) = \prod \limits_{t \in V} \cm{Y}_t(\cm{F})$. Then for all $V \in W$ and $\bar x \in \cm{Y}_V(\cm{F})$, the limits
\begin{equation}
\label{q-Gibbs}
g_V^{\bar x} (x) = \lim \limits_{n \to \infty} \frac{P_{\Lambda_n}(x \bar x_{\Lambda_n \backslash V})}{P_{\Lambda_n \backslash V} (\bar x_{\Lambda_n \backslash V})}, \qquad x \in X^V,
\end{equation}
are positive and the convergence with respect to $\bar x \in \cm{Y}_V(\cm{F})$ is uniform. Moreover, conditional distributions $g_V^{\bar x}$ are quasilocal on sets $\cm{Y}_V(\cm{F})$, $V \in W$.
\end{prop}

\begin{remark}
During the discussion of the first version of the paper, S. Dachian suggested the following equivalent definition of a Gibbs random field.

A random field $P = \{P_V, V \in W\}$ will be called a \emph{Gibbs random field} if for all $t \in \mathbb{Z}^d$, there exists a set of boundary conditions $\cm{Y}_t \subset X^{t^c}$ such that $P(\cm{Y}_t) = 1$ and for any $\bar x \in \cm{Y}_t$, the limits of finite-conditional distributions
$$
\lim \limits_{\Lambda \uparrow t^c} \frac{P_{\Lambda}(x \bar x_\Lambda)}{P_\Lambda (\bar x_\Lambda)} = \lim \limits_{\Lambda \uparrow t^c} g_t^{\bar x_\Lambda}(x) = g_t^{\bar x} (x), \qquad x \in X^t,
$$
exist, are positive and convergence with respect to $\bar x \in \cm{Y}_t$ is uniform.

This definition also can be taken as a basis for the construction of the theory we are developing.  $\quad \bigtriangleup$
\end{remark}

\begin{remark}
In contrast to $\Pr$--definition, in the definition proposed in~\cite{DN09}, it was required that uniform convergence in~\eqref{q1-Gibbs} takes place for all boundary conditions $\bar x \in X^{t^c}$, $t \in \mathbb{Z}^d$. Moreover, in~\cite{DN09}, it was required that uniform convergence in~\eqref{q1-Gibbs} takes place with respect to any increasing sequence $\{\Lambda_n\}_{n \ge 1}$. $\quad \bigtriangleup$
\end{remark}

\begin{remark}
Sullivan~\cite{Sullivan} considered random fields for which limits~\eqref{q-Gibbs} exist for all boundary conditions, are positive and the convergence with respect to all $\bar x$ is uniform. He called such random fields positive almost--Markovian and was looking for corresponding potential. Such potential exists and satisfies a condition which was later called a relatively uniform convergence condition (see, for example,~\cite{LeNy}). Note that the definition of a Gibbs random field suggested in~\cite{DN09}, being similar, was given in terms of one-point conditional distributions only. $\quad \bigtriangleup$
\end{remark}

\begin{remark}
Any random field which is Gibbsian in the sense of the DLR--definition is a Gibbs random field (in the sense of the $\Pr$--definition). Indeed, according to the DLR--definition, a Gibbs random field with uniformly convergent potential $\Phi$ is a random field having a version of its conditional distribution (in Kolmogorov's sense) that almost everywhere coincides with the specification $Q^{\Phi}$. Following our approach, one can formulate the DLR--definition in the next equivalent form: the random field $P =\{P_V, V \in W\}$ is a Gibbs random field corresponding to the uniformly convergent potential $\Phi$ if for all $V \in W$ and $\bar x \in \cm{X}_V(\cm{F})$,
$$
g_V^{\bar x}(x) = \frac{\exp \{-H_V^{\bar x}(x)\}}{\sum \limits_{z \in X^V} \exp \{-H_V^{\bar x}(z)\}}, \qquad x \in X^V,
$$
where $\cm{F} = \{\cm{B}^{\Lambda_n}\}_{n \ge 1}$ is some filtration and $H^\Phi = \{H_V^{\bar x}, \bar x \in X^{V^c}, V \in W\}$ is the Hamiltonian corresponding to potential $\Phi$. The uniform convergence of the potential $\Phi$ implies quasilocality of $g_V^{\bar x}$ on $\cm{X}_V(\cm{F})$ for all $V \in W$, and, particularly, for one-point sets. Hence, by Criterion~\ref{crit-Gibbs-1}, any Gibbs random field corresponding to a uniformly convergent potential $\Phi$ is a Gibbs random field.
\end{remark}

\begin{remark}
There is a definition of a Gibbs random field that is based on the notion of specification and does not use the concept of potential: a random field $P$ is said to be Gibbsian if there is a positive quasilocal specification $Q$ that is a version of the conditional distribution of $P$ (see, for example,~\cite{MaesKo}). Such a definition became possible after the works of Kozlov~\cite{Kozlov} and Sullivan~\cite{Sullivan}, which showed that for a positive quasilocal specification $Q$ there is always a convergent (absolutely uniformly convergent and relatively uniformly convergent, respectively) potential $\Phi$ such that $Q = Q^{\Phi}$. Kozlov also proposed a criterion (see Theorem 2 in~\cite{Kozlov}) for a given random field $P$ to have a version of its conditional distribution that is a positive quasilocal specification. In contrast to it, our Criterion~\ref{crit-Gibbs-1} establishes the Gibbsianness of a random field by imposing conditions directly on its system of conditional distributions. $\quad \bigtriangleup$
\end{remark}

\begin{remark}
$\Pr$--definition, probably, will be useful when considering issues related to Dobrushin's program on Gibbsianity restoration for non-Gibbsian random fields (see, for example,~\cite{vanEnterKo, LeNy}). However, these questions are beyond the scope of our considerations. The $\Pr$--definition is sufficient to achieve the goals set in this paper. $\quad \bigtriangleup$
\end{remark}

\subsection{Representation theorem}

Let $P$ be a random field having a positive one-point conditional distributions $G_1(P,\cm{F})$ with respect to some filtration $\cm{F}$. According to Theorem~\ref{prop-QP-Gibbs}, the system $G_1(P,\cm{F})$ admits a Gibbsian representation in terms of the system $\Delta_1(P,\cm{F}) = \{\Delta_t^{\bar x}, \bar x \in \cm{X}_t(\cm{F}), t \in \mathbb{Z}^d\}$ of one-point transition energies of $P$ with respect to $\cm{F}$.  Since $\Delta_1(P,\cm{F})$ has an autonomously defined analogue --- one-point transition energy field $\boldsymbol{\Delta}_1 = \{ \delta_t^{\bar x}, \bar x \in X^{t^c}, t \in \mathbb{Z}^d \}$, it is possible to consider the inverse problem: under what conditions on $\boldsymbol{\Delta}_1$ there exists a random field $P$ such that the elements of $\boldsymbol{\Delta}_1(P,\cm{F})$ coincide with the corresponding elements of $\boldsymbol{\Delta}_1$, that is, for all $t \in \mathbb{Z}^d$ and $\bar x \in \cm{X}_t(\cm{F})$, it holds
$$
\Delta_t^{\bar x}(x,u) = \delta_t^{\bar x}(x,u), \qquad x,u \in X^t.
$$

The following theorem gives a solution to the mentioned inverse problem.

\begin{T}
\label{th-exist-Gibbs-delta}
Let $\boldsymbol{\Delta}_1 = \{ \delta_t^{\bar x}, \bar x \in X^{t^c}, t \in \mathbb{Z}^d \}$ be a quasilocal one-point transition energy field. Then there exists a random field $P$ corresponding to $\boldsymbol{\Delta}_1$ and $P$ is a Gibbs random field.
\end{T}
\begin{proof}
Any quasilocal one-point transition energy field $\boldsymbol{\Delta}_1$ uniquely determines a quasilocal 1--specification $Q_1$ (see Theorem 2.3 in~\cite{DN19} and a comment after it), which, in its turn, uniquely determines positive quasilocal specification $Q = \{q_V^{\bar x}, \bar x \in X^{V^c}, V \in W\}$ (see Theorem 1.1 in~\cite{DN19}).

Further, according to Dobrushin's existence theorem (see Theorem 1 in~\cite{Dobr}), there exists a random field $P$ whose conditional distributions (in Kolmogorov's sense) $P$-a.e. coincide with the corresponding elements of $Q$. Thus, for any filtration $\cm{F} = \{\cm{B}^{\Lambda_n}\}_{n \ge 1}$, there are sets of admissible boundary conditions $\cm{Y}_t(\cm{F}) \subset \cm{X}_t(\cm{F})$, $P(\cm{Y}_t(\cm{F}))=1$, $t \in \mathbb{Z}^d$, such that for all $t \in \mathbb{Z}^d$ and $\bar x \in \cm{Y}_t(\cm{F})$, it holds
$$
g_t^{\bar x}(x) = q_t^{\bar x}(x), \qquad x \in X^t.
$$
From here it follows that for all $t \in \mathbb{Z}^d$ and $x \in X^t$, $\bar x \in \cm{X}_t(\cm{F})$, we have
$$
g_t^{\bar x_{\Lambda_n}}(x) = \frac{1}{P_{\Lambda_n}({\bar x_{\Lambda_n}})} \int \limits_{\{\bar y \in \cm{Y}_t(\cm{F}): \bar y_{\Lambda_n} = \bar x_{\Lambda_n} \}} g_t^{\bar y}(x) P_{t^c}(d \bar y) = \frac{1}{P_{\Lambda_n}({\bar x_{\Lambda_n}})} \int \limits_{\{\bar y \in \cm{X}_t(\cm{F}): \bar y_{\Lambda_n} = \bar x_{\Lambda_n} \}} q_t^{\bar y}(x) P_{t^c}(d \bar y),
$$
and proceeding as in the proof of Theorem~\ref{GibbsProperty}, we can show that finite-conditional distributions $g_t^{\bar x_{\Lambda_n}}$ converge to $q_t^{\bar x}$ for all $\bar x \in \cm{X}_t(\cm{F})$. Therefore, for all $t \in \mathbb{Z}^d$ and $\bar x \in \cm{X}_t(\cm{F})$, we have
$$
\Delta_t^{\bar x}(x,u) = \ln \frac{g_t^{\bar x}(x)}{g_t^{\bar x}(u)} = \ln \frac{q_t^{\bar x}(x)}{q_t^{\bar x}(u)} = \delta_t^{\bar x}(x,u), \qquad x,u \in X^t,
$$
and hence, random field $P$ corresponds to $\boldsymbol{\Delta}_1$.

Further, for all $\bar x \in \cm{X}_t(\cm{F})$, the limits
$$
g_t^{\bar x}(x) = \lim \limits_{n \to \infty} \dfrac{P_{\Lambda_n} (x \bar x_{\Lambda_n \backslash t})}{ P_{\Lambda_n \backslash t} (\bar x_{\Lambda_n \backslash t}) }, \qquad x \in X^t,
$$
are positive and quasilocal on $\cm{X}_t(\cm{F})$ as functions on $\bar x$. Hence, the convergence with respect to $\bar x \in \cm{X}_t(\cm{F})$ is uniform (see the proof of Theorem~\ref{GibbsProperty}), and therefore, $P$ is a Gibbs random field.
\end{proof}

\begin{remark}
The Gibbs random field, whose existence is established by Theorem~\ref{th-exist-Gibbs-delta}, is Gibbsian in the stronger sense of the definition of works~\cite{DN09,DN19}. Indeed, let $P$ be a random field corresponding to a quasilocal 1--specification~$Q_1$ constructed by $\boldsymbol{\Delta}_1$. Then $Q_1$ is a version of conditional distribution (in Kolmogorov's sense) of the random field $P$. Thus, $P$ has a quasilocal strictly positive 1--specification, and according to Criterion 3.2 in~\cite{DN09}, it is Gibbsian in the sense of the definition of works~\cite{DN09, DN19}. $\quad \bigtriangleup$
\end{remark}

Theorem~\ref{th-exist-Gibbs-delta} is a representation theorem for Gibbs random fields in terms of transition energies. According to it, to construct Gibbs random fields, it is necessary to have a procedure for specifying quasilocal transition energy fields $\boldsymbol{\Delta}_1$.

One of the main ways of such construction is based on the use of potential. Let $\Phi = \{\Phi_V, V \in W\}$ be a uniformly convergent potential and let $H_1^\Phi = \left\{H_t^{\bar x}, \bar x \in X^{t^c}, t \in \mathbb{Z}^d \right\}$ be the corresponding one-point Hamltonian. As is easy to see, the set $\boldsymbol{\Delta}_1^\Phi = \{ \delta_t^{\bar x}, \bar x \in X^{t^c}, t \in \mathbb{Z}^d \}$ of functions
\begin{equation}
\label{Delta-Phi}
\delta_t^{\bar x}(x,u) = H_t^{\bar x}(u) - H_t^{\bar x}(x) = \Phi_t (u) - \Phi_t(x) + \sum \limits_{J \subset W(t^c)} \left(\Phi_{t \cup J}(u \bar x_J) - \Phi_{t \cup J}(x \bar x_J)\right), \qquad x,u \in X^t,
\end{equation}
forms a quasilocal one-point transition energy field. Hence, there exists a Gibbs random field $P$ corresponding to $\boldsymbol{\Delta}_1^\Phi$. The quasilocal transition energy field $\boldsymbol{\Delta}_1^\Phi$ can also be obtained under less restrictive conditions on the potential, namely, it is sufficient to require uniform convergence of the sums in~\eqref{Delta-Phi}. Such potentials were introduced by Sullivan~\cite{Sullivan} and lately was called relatively uniformly convergent (see, for example, Remark 3.58 in~\cite{LeNy}).

Representation theorem is also applicable in the cases when potential is unknown or has a complicated form. For example, to construct a quasilocal transition energy field, one can use the approach proposed in~\cite{DN19} for the discrete Widom--Rowlinson model.

In general, one can proceed as follows. Consider a system $\{\mu_V, V \in W\}$ where $\mu_V$ is a positive measure on $X^V$, $V \in W$. For all $t \in \mathbb{Z}^d$ and $\Lambda \in W(t^c)$, define
$$
\delta_t^{\bar x_\Lambda}(x,u) = \ln \frac{\mu_{t \cup \Lambda}(x \bar x_\Lambda)}{\mu_{t \cup \Lambda}(u \bar x_\Lambda)}, \qquad x,u \in X^t, \bar x \in X^{t^c}.
$$
If for each $\bar x \in X^{t^c}$, there exist uniform (with respect to $\bar x$) limits
$$
\delta_t^{\bar x} (x,u) = \lim \limits_{\Lambda \uparrow t^c} \delta_t^{\bar x_\Lambda}(x,u), \qquad x,u \in X^t,
$$
then the system $\boldsymbol{\Delta}_1^\mu = \{ \delta_t^{\bar x}, \bar x \in X^{t^c}, t \in \mathbb{Z}^d \}$ forms a quasilocal one-point transition energy field.

\subsection{Frames of the theory of Gibbs random fields}

In this section, we present in very general terms the outlines of the theory of Gibbs random fields based on their $\Pr$-definition.

When introducing any class of random processes (fields) into consideration, the following scheme is usually followed. After giving the motivation (springs of action) for its introduction, a definition of the process has to be given in terms of its intrinsic (probabilistic) characteristics. Then, the question of the width of the introduced class of processes is considered, and their general properties are established. Further, a representation theorem is proved that solves the problem of constructing specific processes from this class using well-known or simple in structure non-random objects. This theorem opens up the possibility of constructing as well as studying various processes of this class.

The motivation for introducing Gibbs random fields with a given potential was based on the following thesis: without making a thermodynamic passage to the limit, one can immediately consider limit objects defined in infinite volumes. It was assumed that it is in such idealized systems that the properties of thermodynamic functions manifest themselves most clearly. In particular, this concerns the leading problem of statistical physics --- a strictly mathematical description of the phenomenon of phase transition in physical systems.

The idea for introducing a class of random fields that are Gibbsian in the sense of $\Pr$--definition was (as already noted) our desire to present the theory of Gibbs random fields following the above scheme. Regarding the width of the class of $\Pr$--Gibbs random fields, this question is easily resolved, since the class of positive Markov random fields is its subset. Indeed, the conditional distributions of a Markov random field as a function of boundary conditions are local functions, and therefore, quasilocal. Their uniform limits are also quasilocal, and this allows us to say that the class of Gibbs random fields is a uniform extension of the class of Markov random fields. The width of the class of Gibbs random fields is also given by the fact that any Gibbs random field corresponding to a uniformly convergent potential (that is, random fields which are Gibbsian in the sense of DLR--definition) is Gibbsian by definition.

On the other hand, simple examples show that not all positive random fields are Gibbsian. However, the conditional distributions of any positive random field have a Gibbs form.

For example, there is a general fact according to which any non-trivial mixture of Bernoulli random fields is not Gibbsian (see, for example,~\cite{vEFS}). Below we give an example of such mixture for which finite-dimensional and conditional distributions can be written out explicitly (see~\cite{DN09, DN19} for details), which, in turn, makes it possible to directly verify its non-Gibbsianness. At the same time, due to Theorem~\ref{prop-QP-Gibbs}, conditional distributions of this random field admit Gibbsian representation, and it is possible to give the explicit form of corresponding Hamiltonian. It should be emphasized that the values of this Hamiltonian depend on infinite configurations, and therefore, cannot be represented as a sum of local interactions.

\noindent \textit{Example 2}. Let $\tau>0$ and consider the random field $P = \{P_V,V \in W\}$ with state space $X=\{0,1\}$ which is the mixture with the density $\tau p^{\tau - 1}$ of Bernoulli random fields $P^{(p)}$, $p \in (0,1)$, that is,
$$
P_V(x) = \int \limits_0^1 p^{\vert x \vert} (1-p)^{\vert V \vert - \vert x \vert} \tau p^{\tau-1} dp,
$$
where $\vert x \vert = \vert \{ t \in V : x_t = 1 \} \vert$, $x \in X^V$, $V \in W$.

Let $\cm{F} = \{\cm{B}^{\Lambda_n}\}_{n \ge 1}$ be a filtration. For all $p \in [ 0,1 ]$, denote by $U_t^p(\cm{F})$ the set of configurations $\bar x \in X^{t^c}$ for which the limit $\lim \limits_{n \to \infty} \vert x_{\Lambda_n} \vert /\vert \Lambda_n \vert $ exists and equals $p$. Since
$$
g_t^{\bar x_\Lambda}(1) = \frac{P_{t \cup \Lambda}(1 \bar x_\Lambda)}{P_\Lambda(\bar x_\Lambda)} = \frac{\vert \bar x_\Lambda \vert + \tau}{\vert \Lambda \vert + \tau + 1}, \qquad \bar x \in X^{\mathbb{Z}^d}, \Lambda \in W(t^c), t \in \mathbb{Z}^d,
$$
it can be shown that for all $t\in \mathbb{Z}^d$,
$$
g_t^{\bar x}(x) = p^x (1-p)^{1-x},  \qquad  x \in X^t, \bar x \in U_t^p(\cm{F}), p \in (0,1),
$$
and
$$
g_t^{\bar x}(x) = 1-x, \quad x \in X^t, \bar x \in U_t^0(\cm{F}), \qquad g_t^{\bar x}(x) = x, \quad x \in X^t, \bar x \in U_t^1(\cm{F}).
$$
On the other hand, it is not difficult to verify that for any $\bar x \in X^{\mathbb{Z}^d} \backslash \bigcup \limits_{p \in [0,1]} U_t^p(\cm{F})$, the limits $\lim \limits_{n \to \infty} g_t^{\bar x_{\Lambda_n \backslash t}}(x)$, $x \in X^t$, $t \in \mathbb{Z}^d$, do not exist. Hence, the set $G_1(P,\cm{F}) = \{g_t^{\bar x}, \bar x \in \cm{X}_t(\cm{F}), t \in \mathbb{Z}^d \}$ with $\cm{X}_t(\cm{F}) = \bigcup \limits_{p \in [0,1]} U_t^p(\cm{F})$ forms a system of conditional distributions of the random field $P$ with respect to filtration~$\cm{F}$.

Now, for each $t \in \mathbb{Z}^d$, consider the set $\cm{X}_t = \bigcap \limits_{\cm{F}} \cm{X}_t(\cm{F})$. It is not difficult to see that $\cm{X}_t$ consists of configurations $\bar x \in X^{t^c}$ for which there exists the limit $\lim \limits_{\Lambda \uparrow t^c} \vert x_\Lambda \vert /\vert \Lambda \vert $. It can be shown that
$$
\cm{X}_t = \{\bar x \in X^{t^c}: \vert \{ s \in t^c: \bar x_s = 1 \} \vert < \infty \} \cup \{\bar x \in X^{t^c}: \vert \{ s \in t^c: \bar x_s = 0 \}\vert < \infty \}.
$$
Therefore, $P(\cm{X}_t) = 0$ for any $t \in \mathbb{Z}^d$, and according to Proposition~\ref{nonGibbs}, $P$ is not a Gibbs random field.

At the same time, allowing infinite values for the Hamiltonian (which does not contradict its axiomatic definition), it can be shown that the elements of $G_1(P,\cm{F})$ admit a Gibbs form with a one-point Hamiltonian $H_1(P,\cm{F}) = \{H_t^{\bar x}, \bar x \in \cm{X}_t(\cm{F}), t \in \mathbb{Z}^d \}$ defined by
$$
H_t^{\bar x}(x) = \left\{ \begin{array}{l}
                              - x \ln p - (1-x) \ln(1-p), \qquad x \in X^t, \bar x \in U_t^p(\cm{F}), p \in (0,1),\\
                              \\
                              0, \qquad \quad x = 0, \bar x \in U_t^0(\cm{F}) \text{ or } x = 1, \bar x \in U_t^1(\cm{F}),\\
                              \\
                              + \infty, \qquad x = 1, \bar x \in U_t^0(\cm{F})  \text{ or } x = 0, \bar x \in U_t^1(\cm{F}).
                  \end{array} \right.
$$
\qed\\

Thus, the set of Gibbs random fields is a proper subset of the set of all positive random fields. At the same time, it can be shown that the class of Gibbs random fields is dense in the set of all positive random fields with respect to the weak topology (see~\cite{DN09, DN19}).

Now about the general properties of $\Pr$--Gibbs random fields. These include criteria for Gibbsianness of random fields and properties of weak dependence (ergodicity, regularity, uniform strong mixing, semi-invariant mixing conditions, etc.). Similar properties also include almost the entire spectrum of limit theorems of probability theory (central and local limit theorems, estimates of the rate of convergence and asymptotic expansions in them, the law of the iterated logarithm, etc.), and, especially the important functional limit theorem (see, for example,~\cite{N91}). In these problems, one can also use the martingale approach, which allows proving limit theorems under weaker assumptions, for example, based on the ergodicity property (see~\cite{NP92, N97}).

The application of limit theorems in statistical physics is fundamentally important when considering the problem of equivalence of ensembles and the theory of fluctuations of physical quantities. From a physical point of view, it is clear that various ensembles of statistical physics in the thermodynamic limit should equally describe the corresponding physical systems. The mathematical justification for this fact is the problem of equivalence of ensembles. The main tool for proving the equivalence of ensembles for any physical system is local limit theorems for the number of particles and energy (see~\cite{MKh}).

Macroscopic physical quantities are usually taken equal to their average values. In this case, fluctuations occur: deviations of the real values of a physical quantity from its average value, as a result of which the problem of estimating the probability of these deviations arises. The answer to this question is given by the central limit theorem and its refinements.

In a number of cases, some important properties of thermodynamic functions have both physical and probabilistic interpretations. For example, the property of convexity of free energy (or pressure) is equivalent to the linear growth (with respect to volume) of the dispersion of the Hamiltonian (see~\cite{DN}). Another typical example is Dobrushin's probabilistic criterion for the presence of a phase transition in a physical system. According to this criterion, the non-uniqueness of a Gibbs random field with a given potential can be interpreted as the presence of a phase transition in the physical system (see~\cite{DUniq}).

All these general results are quite applicable to specific Gibbs random fields. The latter are built on the basis of the representation theorem.

\section{Concluding remarks}

Mathematical statistical physics can be considered from two positions: physical and probabilistic. The DLR--definition, reflecting the physical point of view, allows within its framework to build and study mathematical models of physical systems based on the notion of potential. This definition is not simple even for specialists in probability theory since for its perception, it requires a comprehension of such physical concepts as potential, Hamiltonian, and Gibbs specification. The $\Pr$--definition, reflecting the probabilistic point of view, has a simpler formulation, based solely on the properties of finite-dimensional distributions of a random field without appeal to any physical notions. Gibbs random field in the sense of the DLR--definition is such in the sense of the $\Pr$--definition. The $\Pr$--definition reveals the probabilistic meaning of the class of Gibbs random fields as a uniform extension of the class of Markov random fields.

The symbiosis of these two definitions of Gibbs random fields allows us to present a consistent version of the theory of Gibbs random fields if we take the $\Pr$--definition as the main one, and consider the existence theorem within the framework of the DLR--definition as a representation theorem.\\

\noindent \textbf{Acknowledgments.} The authors are grateful to Serguei Dachian, University of Lille, for his careful reading of the first version of the paper and very significant comments. We are also grateful to Aernout van Enter for his useful remarks. The work was supported by the Science Committee of the Republic of Armenia in the frames of the research project 21AG-1A045.

\end{document}